# Survival exponents for fractional Brownian motion with multivariate time

## G. Molchan

Institute of Earthquake Prediction Theory and Mathematical Geophysics,

Russian Academy of Science, 84/32 Profsoyuznaya st.,

117997, Moscow, Russian Federation

E-mail address: molchan@mitp.ru

*Abstract*. Fractional Brownian motion, H-FBM, of index $H \in (0,1)$ with d-dimensional time is considered in a spherical domain that contains 0 at its boundary. The main result: the log-asymptotics of probability that H-FBM does not exceed a fixed positive level is (H-d)logT(1+o(1)), where T>>1 is radius of the domain.

## 1. Introduction.

Fractional Brownian motion, H-FBM of index $H \in (0,1)$ with multivariate time $t \in R^d$ is a centered Gaussian random process $w_H(t)$ with correlation function

$$Ew_H(t)w(s) = 1/2(|t|^{2H} + |s|^{2H} - |t-s|^{2H}).$$

H-FBM is H-self-similar (H-ss), isotropic, and has stationary increments (si), i.e.,

$\{w_H(\lambda Ut + t_0) - w_H(t_0)\} \doteq \{\lambda^H w_H(t)\}$ holds in the sense of the equality of finite-dimensional distributions for any fixed $t_0$, $\lambda > 0$, and rotation around 0, $U$.

The one-sided exit problem for a random process $\xi(t)$ and its characteristics, the so-called survival exponents:

$$\theta_\xi = \lim_{T \to \infty} -\log P(\xi(t) < 1, t \in \Delta_T)/\psi(T) \qquad (1)$$

are the subject of intensive analysis in applications. Here $\Delta_T$ is an increasing sequence of domains of size $T$ and $\psi(\cdot)$ is a suitable slowly varying function, typically, $\psi(T) = \log T$ for ss-processes. The greatest progress in this area has been achieved for processes with one-dimensional time. (See surveys by Bray et al (2013) in the physics literature and by Aurzada and Simon (2015) in the mathematical one).

2H-FBM was one of the first non-trivial examples of non-Markovian processes for which the survival exponents have been found exactly (Molchan,1999). Namely, for $\psi(T) = \log T$, the survival exponents are

$$\theta_{w_H} = 1 - H \ , \ \Delta_T = (0,T) \quad \text{and} \quad \theta_{w_H} = d \ , \ \Delta_T = (-T,T)^d \ . \tag{2}$$

Recently, Aurzada et al (2016) considerably refined the asymptotics of probability

$$p_T = P(w_H(t) < 1, t \in \Delta_T) \ , \quad \Delta_T = (0,T) \ , \tag{3}$$

and showed that the exponent $\theta = 1 - H$ is universal for a broad class of H-ss processes with stationary increments. The ideas of this work have proved useful in analysis of the conjunction that $\theta_{w_H} = d - kH$ for $w_H(t)$ in $\Delta_T = [0,T]^k \times [-T,T]^{d-k}$ (Molchan , 2012).

The case $k = 0$ corresponds to the right part of (2). The case $k = 1$ is supported by the result which we discuss below: $\theta_{w_H} = d - H$ for fractional Brownian motion in $\Delta_T = T\Delta_1$ where $\Delta_1$ is a unit ball that contains 0 at its boundary.

To estimate $\theta_{w_H}$ we modify the famous result by Aurzada et al (2016), which says that, for a broad class of si-processes: $\xi(t), \xi(0) = 0$ with discrete time, $t \in Z^1$,

$$|\Delta_T| P(\xi(t) < 1, t \in \Delta_T \setminus \{0\}) \propto E \max(\xi(t), t \in \Delta_T) \ , \tag{4}$$

where $\Delta_T = [0,T]$ and $|\Delta_T| = T$ .

For H-ss processes with continuous time, the right part of (4) is proportional to $T^H$, and therefore the exponent for (3) is $1 - H$.

However, the result by Aurzada et al (2016) essentially uses the 1-D nature of time. Considering $|\Delta_T|$ as the volume of $\Delta_T$, relation (4) is found to be in formal agreement with the conjunction for $k = 1$ but not for $k > 1$; in addition, (4) is very crude for $k = 0$ (see (3)). This means that the analysis of the cases $d > 1, k > 1$ needs in additional ideas.

## 2. The lower bound.

**Proposition 1.** Let $\xi(t), \xi(0) = 0, t \in R^d$ be a centered isotropic random process with stationary increments. Then

$$P(\xi(t) < -1, t \in \Delta_T, |t| > 1) < cT^{-d} E \max(\xi(t), t \in \Delta_T) \ , \tag{5}$$

where $\Delta_T = T\Delta_1$ is a ball of radius $T$ that contains 0 at its boundary.

***Consequence***. If $\xi(t)$ is fractional Brownian motion of index $H \in (0,1)$ in $\Delta_T$, then the survival exponent has the lower bound $\theta_{w_H} \geq d - H$.

***Remark***. Proposition 1 holds for $\Delta_T = [0,T] \times [-T,T]^{d-1}$ as well.

**Proof.** Let $U_T = \{x_{k,\alpha}, \alpha = 1,2...n_k; k = 1,2...\}$ be a subset of ball $B_T$ of radius $T$ in $R^d$; $U_T$ consists of $N_T$ points such that

$$|x_{k,\alpha}| = r_k, \quad |x_{k,\alpha} - x_{m,\beta}| > 1, \quad N_T > CT^d, \quad 1 < r_k < r_{k+1} \leq T. \tag{6}$$

Consider the following increasing sequence of subsets of $U_T$:

$$U_{k+1,\alpha} = U_k \cup \bigcup_{\beta=1}^{\alpha} x_{k+1,\beta}, \quad U_k = \{x_{i,\gamma} : |x_{i,\gamma}| \leq r_k\}.$$

Fix $\Delta_T = \{t : |t + Te| = T\}$, where $e = (0,0....,1)$. Let $O_{k,\alpha}$ be a rotation transferring $x_{k,\alpha}$ in $\tilde{x}_{k,\alpha} = r_k e$. Setting $\tilde{U}_{k,\alpha} = O_{k,\alpha} U_{k,\alpha}$, one has

$$(\tilde{U}_{k,\alpha} - \tilde{x}_{k,\alpha}) \setminus \{0\} \subset \Delta_T \setminus B_1, \quad (k,\alpha) \neq (1,1) \tag{7}$$

Therefore, using the notation $M(A) = \sup(\xi(t), t \in A)$, we get

$$p_T(-1) := P(\xi(t) < -1, t \in \Delta_T \setminus B_1) \leq P(M((\tilde{U}_{k,\alpha} - \tilde{x}_{k,\alpha}) \setminus \{0\}) < -1). \tag{8}$$

By the si-property of $\xi(t)$, we can continue

$$= P(M(\tilde{U}_{k,\alpha} \setminus \tilde{x}_{k,\alpha}) - \xi(\tilde{x}_{k,\alpha}) < -1) = P(M(U_{k,\alpha-1}) + 1 < \xi(x_{k,\alpha})). \tag{9}$$

The last equality holds because $\xi(t)$ is rotation invariant.

The event $\{M(U_{k,\alpha-1}) + 1 < \xi(x_{k,\alpha})\}$ is measurable relative to the sequence

$$\xi(x_{1,1}),...\xi(x_{1,n_1});........\xi(x_{k,1}),...\xi(x_{k,n_k})...... \tag{10}$$

and means that $\xi(x_{k,\alpha})$ is a record which exceeds the previous one by at least 1. Let $\nu_T$ be the number of such records in (10). Then, by (8,9),

$$(N_T - 1) p_T(-1) \leq \sum_{k,\alpha} P(M(U_{k,\alpha}) + 1 < \xi(x_{k,\alpha+1})) = E\nu_T \leq E(M(U_T) - \xi(x_{1,1})),$$

where $U_{1,1} = \{x_{1,1}\}, \quad U_{k,n_k} = U_{k+1}, \quad x_{k,n_k+1} = x_{k+1,1}, \quad (k,\alpha) \neq (1,1).$

Finally, by (6),

$$p_T(-1) \leq E(M(U_T))/(N_T - 1) < cT^{-d} E(\sup \xi(t), t \in \Delta_T). \tag{11}$$

Suppose that $\xi(t)$ is fractional Brownian motion of index $H \in (0,1)$ in $\Delta_T$. By the standard procedure, we can compare $p_T(-1)$ with

$$p_T(1) = P(w_H(t) < 1, t \in (\Delta_T \setminus B_1)). \tag{12}$$

For this purpose we can find a continuous function $\varphi_T(t)$ such that

$$\varphi_T(t) = 1, |t| > 1, \ \|\varphi_T\|_{H,T} < const, \tag{13}$$

where $\|\cdot\|_{H,T}$ is the norm of the Hilbert space $H_H(\Delta_T)$ with the reproducing kernel $Ew_H(t)w_h(s), (t,s) \in \Delta_T \times \Delta_T$ (see for this fact Molchan(1999) or Appendix). Then

$$p_T(-1) = P(w_H(t) + 2\varphi_T(t) < 1, t \in (\Delta_T \setminus B_1)).$$

According to (Aurzada&Dereich, 2010),

$$\left| \sqrt{\ln 1/p_T(1)} - \sqrt{\ln 1/p_T(-1)} \right| \leq \|2\varphi_T\|_{H,T} / \sqrt{2}. \tag{14}$$

From the self-similarity of H-FBM and (11) one has

$$p_T(-1) \leq cT^{-(d-H)} EM_{w_H}(\Delta_1). \tag{15}$$

Combining (13-15), one has

$$[\ln 1/P(w_H(t) < 1, t \in \Delta_T)]^{1/2} / \sqrt{\ln T} \geq \sqrt{d-H} + O(1/\sqrt{\ln T}), \tag{16}$$

i.e., $\theta_{w_H} \geq d - H$.

## 3. The upper bound.

Below we use notation $M(A) = \sup(w_H(t), t \in A)$ and $|A| = \#\{t : t \in A\}$.

**Proposition 2.** Let $w_H(t)$ be H-FBM in $\Delta_T = T\Delta_1 \subset R^d$ where $\Delta_1$ the is a bounded domain and $0 \in \Delta_1$. Consider a finite 1-net of $\Delta_T$, i.e. a subset $U_T = \{x_k, k=1,...,N_T\} \subset \Delta_T$, $\{0\} \notin U_T$ such that

$$N_T \propto T^d \text{ and } \Delta_T \subset \cup_1^{N_T} B_1(x_r),$$

where $B_1(x)$ is unite ball with center $x$. Then for $T > T_0$ and $0 < q < 1$

$$P(M(\Delta_T) < 2c_H\sqrt{d \ln T}) > qP(M(U_T) < 0), \tag{17}$$

where $c_H = 1 + (2+\sqrt{2})\int_1^\infty 2^{-Hv^2} dv$ is the Fernique constant.

In addition,

$$EM(U_T) = EM(\Delta_T)(1+o(1)) = T^H EM(\Delta_1)(1+o(1)) \ , T \to \infty \ . \qquad (18)$$

**Proof.** One has

$$P(M(U_T) < 0) < P(M(U_T) < 0, A_T) + P(A_T^c) \ , \qquad (19)$$

where

$$A_T = \{\max_k \max_t (w_H(t) - w(x_k), t \in B_1(x_k)) < b_T\} .$$

We can continue the previous inequality

$$< P(M(\Delta_T) < b_T) + \sum_k P(\max(w_H(t) - w_H(x_k), t \in B_1(x_k)) > b_T)$$

$$< P(M(\Delta_T) < b_T) + N_T P(M(B_1) > b_T) := p_{1,T} + p_{2,T} . \qquad (20)$$

Applying the Fernique (1975) result to $w_H(t)$, we have

$$P(M(B_1) > r_T c_H) < c_d \int_{r_T}^\infty e^{-u^2/2} du \ , \quad r_T > (1+4d)^{1/2} \ . \qquad (21)$$

From here, setting $b_T = \sqrt{4d \ln T} c_H$, one has

$$p_{2,T} < CT^d \cdot T^{-2d}/\sqrt{\ln T} = CT^{-d}/\sqrt{\ln T} \ . \qquad (22)$$

To show $p_{2,T} = o(p_{1,T})$, note that $\Delta_T \subset B_{TD}$, where D is diameter of $\Delta_1$. Therefore

$$p_{1,T} = P(M(\Delta_T) < b_T) > P(M(B_{TD}) < b_T) = P(M(B_{T'}) < 1), \quad T' = TD/b_T^{1/H} . \qquad (23)$$

By Molchan (1999), $P(M(B_T) < 1) = T^{-d}(1+o(1))$. Due to (22),(23), we have

$$p_{2,T}/p_{1,T} = O(\ln T)^{-(1+d/H)/2}) = o(1). \qquad (24)$$

Relations (19, 20) and (24) imply (17).

To prove relation (18), note that

$$M(\Delta_T) < M(U_T) + \max_k \max_t (w_H(t) - w(x_k), t \in B_1(x_k)) := M(U_T) + \delta_T . \qquad (25)$$

As above, using the event $A_T = \{\max_k \max_t (w_H(t) - w(x_k), t \in B_1(x_k)) < b_T\}$, one has

$$E\delta_T < b_T + E\delta_T 1_{A_T} < b_T + N_T EM(B_1)[M(B_1). > b_T], \qquad (26)$$

where $b_T = \sqrt{4d \ln T} c_H$ and $N_T < CT^d$. Therefore, the 2-d term in (26) is $o(1)$.

By (25), (26), we obtain (18), because

$$EM(U_T) > EM(\Delta_T) - E\delta_T > EM(\Delta_T) - c\sqrt{\ln T} + o(1) = T^H EM(\Delta_1) - c\sqrt{\ln T} + o(1) .$$

**Proposition 3.** Let $w_H(t), t \in \Delta_T$ be H-FBM, $\Delta_T = T\Delta_1 \subset R^d$, where $\Delta_1$ is a unite ball and $0 \in \Delta_1$. Then

$$P(M(\Delta_T) < 1) > c(T\sqrt{\ln T})^{-(d-H)},$$

i.e. the survival exponent for H-FBM in $\Delta_T$ has the upper bound $\theta_{w_H} \leq d - H$.

***Consequence***. Due to Propositions 1, 3, the survival exponent for H-FBM in $\Delta_T$ exists and is equal to $d - H$.

**Proof**. As in proof of Proposition 1, we consider again the subset $U_T$ of ball $B_T \subset R^d$:
$U_T = \{x_{k,\alpha}, \alpha = 1,2...n_k; k = 1,2...\}$, $\{0\} \notin U_T$ In addition to the properties (6), we suppose that the elements of $U_T$ are numerated in such way that

$$x_{k,\alpha+1} \in B_2(x_{k,\alpha}) \text{ and } x_{k+1,1} \in B_2(x_{k,n(k)}). \qquad (27)$$

As before,

$$U_{k+1,\alpha} = U_k \cup \bigcup_{\beta=1}^{\alpha} x_{k+1,\beta} \quad , \quad U_k = \{x_{i,\gamma} : |x_{i,\gamma}| \leq r_k\} := U_{k,0} ;$$

$\Delta_T = \{t : |t + Te| = T\}$, where $e = (0,0....,1)$; $O_{k,\alpha}$ is a rotation transferring $x_{k,\alpha}$ in $\widetilde{x}_{k,\alpha} = r_k e$. Setting $\widetilde{U}_{k,\alpha} = O_{k,\alpha} U_{k,\alpha}$, one has

$$(\widetilde{U}_{k+1,\alpha} - \widetilde{x}_{k+1,\alpha}) \setminus \{0\} \subset \Delta_{k+1} \setminus B_1 .$$

Due to (27), $(\widetilde{U}_{k+1,\alpha} - \widetilde{x}_{k+1,\alpha})$ is 2-net in $\Delta_{k+1}$. Therefore, by (17), for $k > T_0$

$$P(M(\Delta_k) < 2c_H \sqrt{d \ln k}) > qP(M(\widetilde{U}_{k,\alpha} - \widetilde{x}_{k,\alpha}) \setminus \{0\}) < 0)$$

$$= qP(M(\widetilde{U}_{k,\alpha-1}) - w_H(\widetilde{x}_{k,\alpha})) < 0) = qP(M(U_{k,\alpha-1}) < w_H(x_{k,\alpha})).$$

As a result,

$$\sum_{k=K}^{K'} n_k P(M(\Delta_k) < 2c_H \sqrt{d \ln k}) > q \sum_{k=K}^{K'} \sum_{\alpha=1}^{n_k} P(M(U_{k,\alpha-1}) < w_H(x_{k,\alpha})) \qquad (28)$$

where $K = [T]$ and $K' = [T']$.

Similarly to the proof of Proposition 1, we conclude, that the right part of (28) is equal to $E\nu(T,T')$, where $\nu(T,T')$ is a number of records in the following sequences:

$$M(U_K), w_H(x_{K+1,1}), ..., w_H(x_{K+1,n(K+1)_1}); ........ w_H(x_{K',1}), ...w_H(x_{K',n(K')}).$$

Let $\delta(T,T')$ be the maximum increment between adjacent elements of the sequence

$$w_H(x_{K,n(K)}), w_H(x_{K+1,1}), ..., w_H(x_{K+1,n(K+1)_1}); ........ w_H(x_{K',1}), ...w_H(x_{K',n(K')}).$$

Then

$$M(U_{K'}) - M(U_K) \le (\nu(T,T')+1)\delta(T,T') \le (\nu(T,T')+1)b_T + R_T, \qquad (29)$$

where

$$R_T = (|U_{K'} \setminus U_K|+1)\delta(T,T')[\delta(T,T') > b_T].$$

Due to (27),

$$ER_T < (|U_{K'} \setminus U_K|+1)^2 \max_{|t|<2} Ew_H(t)[w_H(t) > b_T].$$

Setting $b_T = 2\sqrt{n \ln T}$ and $T' - T = \rho T$ we obtain

$$ER_T < cT^{2d} \cdot T^{-n} = cT^{2d-n}. \qquad (30)$$

By (29),

$$E\nu(T,T') > EM(U_{K'}) - EM(U_K) - b_T - ER_T,$$

where, according to (18),

$$EM(U_K) = K^H EM(\Delta_1)(1 + o(1)).$$

Setting $n = 2d$, one has

$$E\nu(T,T') > c(T^H - \sqrt{\ln T} - T^{2d-n}) = cT^H(1+o(1)). \tag{31}$$

Now we can continue (28) as follows:

$$qE\nu(T',T) < \sum_{k=K}^{K'} n_k P(M(\Delta_k) < 2c_H\sqrt{d \ln k}). \tag{32}$$

Due to ss-propery of H-FBM,

$$P(M(\Delta_k) < c\sqrt{\ln k}) = P(M(\Delta_{\tilde{k}}) < 1) \quad , \quad \tilde{k} = k/(c\sqrt{\ln k})^{1/H},$$

and therefore the probability term decreases as function of k. Hence, (32) implies

$$qE\nu(T',T) < |U_{T'} \setminus U_T| P(M(\Delta_{T'+1}) < 2c_H\sqrt{d \ln(T'+1)}) < CT^d P(M(\Delta_{\tilde{T}}) < 1), \tag{33}$$

where

$$\tilde{T} = T'/(2c_H\sqrt{d \ln T'})^{1/H} \quad \text{or} \quad T' = \tilde{T}(2c_H\sqrt{d \ln \tilde{T}})^{1/H}(1+o(1)).$$

Finally, by (31) and $T' - T = \rho T$,

$$P(M(\Delta_{\tilde{T}}) < 1) > c(\tilde{T}\sqrt{\ln \tilde{T}})^{-(d-H)}.$$

**Appendix**

***Example from Proposition 1.*** Consider H-FBM in domains $\Delta_T = T \cdot \Delta_1, 0 \in \partial\Delta_1$; then a suitable function $\varphi_T(t), t \in \Delta_T$ can be chosen as follows:

$$\varphi_T(t) = f(|t|/(Tk)) - f(|t|),$$

where $f(x), x \in R^1$ is a finite smooth function such that $f(t) = 1$ for $|x| < 1/2$ and $f(t) = 0$ for $|x| > 1$. Here $k$ is the diameter of $\Delta_1$.

This can be seen as follows (Molchan, 1999). Due to the spectral representation of H-FBM, the Hilbert space $H_H(\Delta_T)$ with the reproducing kernel $Ew_H(t)w_h(s), (t,s) \in \Delta_T \times \Delta_T$ (Aronszajn, 1950), is closure of smooth functions $\varphi(t), \varphi(0) = 0$ relative to the norm

$$\|\varphi\|_{H,T} = \inf_{\tilde{\varphi}} \|\tilde{\varphi}\|_H, \quad \|\psi\|_H = c_H \int |\hat{\psi}(\lambda)|^2 |\lambda|^{d+2H} d\lambda,$$

Where $\tilde{\varphi}(t)$ is a finite function such that $\tilde{\varphi}(t) = \varphi(t), t \in \Delta_T$; $\hat{\psi}(\lambda), \lambda \in R^d$ is the Fourier transform of $\psi(t)$. Obviously, we have $\varphi_T(0) = 0$, $\varphi_T(1) = 1$ for $t \in \Delta_T \setminus B_1$, and

$$\|\varphi\|_{H,T} < \|f(|t|/Tk) - f(|t|)\|_H < \|f(|t|/Tk)\|_H + \|f\|_H = ((Tk)^{-H} + 1)\|f\|_H < 2\|f\|_H.$$